\documentclass{amsart}

\usepackage{amssymb,latexsym,amsmath,amsthm,amsfonts}

\newcommand{\R}{\mathbb{R}}
\newcommand{\C}{\mathbb{C}}

\newcommand{\N}{\mathbb{N}}

\newcommand{\NN}{\mathcal{N}}
\newcommand{\PP}{\mathcal{P}}
\newcommand{\K}{\mathcal{K}}
\newcommand{\tr}{\text{ tr}}

\newtheorem{theorem}{Theorem}[section]

\theoremstyle{definition}

\newtheorem{example}[theorem]{Example}

\numberwithin{equation}{section}

\begin{document}

\title[A Reproducing Kernel Condition for Indeterminacy]{A Reproducing Kernel
		Condition for Indeterminacy in the Multidimensional Moment Problem\\
		\date{April 19, 2006}}
\author[R. A. Roybal]{Roger A. Roybal}
	\address{Department of Mathematics\\
		CSU Channel Islands\\
		One University Drive\\
		Camarillo, CA 93012 USA}
	\email{roger.roybal@csuci.edu}

\thanks{The author would like to thank Mihai Putinar for all his advice and support during the preparation of this
	article.}

\thanks{This paper is to appear in the Proceedings of the American Mathematical Society.}



\begin{abstract}
	Using the smallest eigenvalues of Hankel forms associated with a multidimensional moment problem, we establish a
	condition equivalent to the existence of a reproducing kernel.  This result is a multivariate analogue of Berg, Chen,
	and Ismail's 2002 result.  We also present a class of measures for which the existence of a reproducing kernel
	implies indeterminacy.
\end{abstract}

\maketitle

\section{Introduction}

In \cite{BCI}, Berg, Chen, and Ismail find a new condition equivalent to determinacy in the one-dimensional moment
problem.

\begin{theorem}[Berg, Chen, and Ismail, 2002] \label{T:BCI}
	Let $\lambda_N$ be the smallest eigenvalue of the truncated Hankel matrix $H_N$ for the measure $\mu$.  Then
	$\lambda_N \to 0$ as $N \to \infty$ if and only if $\mu$ is determinate.
\end{theorem}

They use the classical fact that a measure is indeterminate if and only if a reproducing kernel exists on the Hilbert
space in which the polynomials are dense.  Additionally, a reproducing kernel exists if and only the sum
\[
	\sum_{k=0}^\infty |P_k(z_0)|^2
\]
converges for some $z_0 \in \C \backslash \R$, where $\{P_k\}$ is the standard set of orthonormal polynomials.  These
results go back to M. Riesz and may be found in \cite{A}.

Let $\mu$ be a positive measure such that $s_n = \int_\R x^n \, d\mu$ is finite for all $n$.  For each $N \in \N$, we
define the $N^\text{th}$ Hankel matrix to be
\[
	H_N = (s_{i+j})_{i,j=0}^N.
\]
Since $\mu$ is positive, this implies that $H_N$ is positive semidefinite for each $N$.  We let $\lambda_N$ be the smallest
eigenvalue of $H_N$ and note that Cauchy's interlace theorem implies that $\lambda_N$ decreases as $N$ increases.  If
$\lambda_N = 0$ for some $N$, then $\lambda_n = 0$ for all $n \geq N$, and the measure is a finite sum of point masses
\cite{A}.
This case implies determinacy since any measure with compact support is automatically determinate.  In the case where $\mu$
has infinite support, then $\lambda_N > 0$ for every $N$.

Berg, Chen, and Ismail show that $\lambda_N \to \gamma > 0$ if and only if a reproducing kernel exists for the polynomials,
which gives the conclusion.  We propose to extend this to a similar question in the multidimensional setting.

Let $\alpha = (\alpha_1, \alpha_2, \ldots, \alpha_d) \in \N_0^d$ be a multi-index, $x = (x_1, x_2, \ldots, x_d) \in \R^d$,
and let $x^\alpha = x_1^{\alpha_1} x_2^{\alpha_2} \cdots x_d^{\alpha_d}$.  We shall use the standard notation of $|\alpha|
= \alpha_1 + \cdots + \alpha_d$.  For a given multisequence $\{ s_\alpha \}$, define a linear functional $L : \C[x] \to \C$
by $L(x^\alpha) = s_\alpha$, and extend linearly.
We say the multisequence $\{s_\alpha\}$ and the associated functional $L$ are positive if
for every polynomial $p \in \C[x]$, we have $L(p \overline{p}) \geq 0$.  When $L$ is positive, we construct a pre-inner
product by $\langle p,q \rangle = L(p \overline{q})$.  If we consider the ideal $\NN = \{ p \in \C[x] : \langle p,p \rangle
=0 \}$ then $\langle \cdot, \cdot \rangle$ acting on $\C[x] / \NN$ becomes a postive definite form, and we complete
this to a Hilbert space $\PP$ in which $\C[x] / \NN$ is dense.  This is the standard GNS construction and may be found in
full detail in \cite{F}.
We will also assume a normalizing condition of $s_{(0,0,\ldots,0)} = 1$ meaning that we deal with probability measures.
This assumption also simplifies many calculations.

We define the Hankel kernel indexed over $\alpha$ associated with the multisequence $\{ s_\alpha \}$ by
\[
	H = ( s_{\alpha + \beta} )_{\alpha, \beta \in \N_0^d},
\]
and for any $N \in \N$ the $N^\text{th}$ truncation of $H$ by
\[
	H_N = ( s_{\alpha + \beta} )_{0 \leq |\alpha|, |\beta| \leq N}.
\]

One significant difference between the multidimensional and one dimensional cases is that in one variable, a positive
measure $\mu$ exists in $\R$ which represents $L$ in the sense that $L(p) = \int_\R p \, d\mu$ if and only if $L$ is
positive, whereas there exist positive $L$ for which there are no such $\mu$ on $\R^d$ if $d \geq 2$, examples of which
were concurrently discovered in \cite{BCJ} and
\cite{Sch}.  In the case where such a $\mu$ exists, its support must lie in the real algebraic variety generated by the
ideal $\NN$.  If we denote by $\lambda_N$ the smallest eigenvalue of the matrix $H_N$, then $\lambda_N \geq 0$ for all $N
\in \N$ and equality holds for some $N$ if and only if $\NN$ is a nontrivial ideal of $\C[x]$.  For the remainder of this
article we will assume that $\NN = (0)$ so that $\lambda_N > 0$ for all $N$.

If we try to bring Berg, Chen, and Ismail's result to the multivariate case we must find a new statement since there exist
indeterminate measures $\mu$ for which $\lambda_N \to 0$ as $N \to \infty$.

To show this we use a theorem of Petersen~\cite{P}.

\begin{theorem}[L. C. Petersen]
	Let $\mu$ be a moment measure on $\R^d$.  Consider the projection $\pi_j: \R^d \to \R$ onto the $j^{th}$ coordinate
	and the projection measures $\mu_j := \pi_{j *}(\mu)$ on $\R$.
	\begin{itemize}
		\item[(a)] If each $\mu_j$ is determinate, then so is $\mu$.
		\item[(b)] If we assume that $\mu = \mu_1 \otimes \mu_2 \otimes \cdots \otimes \mu_d$, then $\mu$ is
			determinate if and only if $\mu_j$ is determinate for $1 \leq j \leq d$.
	\end{itemize}
\end{theorem}

\begin{example}
	Let $\mu_1$ and $\mu_2$ be measures on $\R$ admitting all moments with $\mu_1$ determinate and $\mu_2$ indeterminate.
	Let $\nu$ be a measure on $\R$ distinct from $\mu_2$ but possessing the same moments.  Define $\mu := \mu_1 \otimes
	\mu_2$.  Clearly $\mu$ is indeterminate since $\mu_1 \otimes \nu$ gives the same moments as $\mu$.  We now show that
	for $\mu$, $\lambda_N \to 0$ as $N \to \infty$.

	Consider the matrix $J_N = (s_{(m + n,0)})_{0 \leq m,n \leq N}$, and let $\eta_N$ be its smallest eigenvalue.  Then
	$J_N$ is a principal submatrix of the matrix $H_N$, and by their self-adjointness and Cauchy's
	interlace theorem we see that $0 \leq \lambda_N \leq \eta_N$.  Note that $J_N$ is the $N^\text{th}$ Hankel matrix
	associated with the measure $\mu_1$ which is determinate in $\R$, whence $\eta_N \to 0$ as $N \to \infty$ by
	\cite{BCI}.  Thus $\mu$ is an indeterminate measure such that $\lambda_N \to 0$.
\end{example}

A reproducing kernel for $\PP$, the completion of the polynomials on $\C^d$, is a function $K : \C^d \times \C^d \to \C$
which has the reproducing property, i.e.\ so that for any polynomial $p(x)$,
\[
	p(y) = \langle p(x), K(x,y) \rangle_\PP.
\]
We construct a complete system of orthonormal polynomials $\{ P_\alpha \}_{\alpha \in \N_0^d}$ which serve as an
orthonormal basis of $\PP$.  Typically one works with orthonormal polynomials in which $\deg(P_\alpha)$ is an
increasing with $|\alpha|$, and here we assume that $\deg(P_\alpha) = |\alpha|$.  This may be achieved by ordering
$\{ x^\alpha \}_{\alpha \in \N_0^d}$ in a graded lexicographical order then using the Gram-Schmidt procedure to obtain an
orthonormal system of polynomials.  The construction of orthonormal polynomials from the monomials may be found in
\cite{DX}.

Using this basis of orthonormal polynomials, we find
\[
	K(x,y) = \sum_{\alpha \in \N_0^d} P_\alpha(x) \overline{P_\alpha (y)},
\]
and for a fixed $y \in \C^d$, this is a function in $\PP$ provided
\begin{equation} \label{E:sum}
	\sum_{\alpha \in \N_0^d} |P_\alpha (y)|^2
\end{equation}
is finite.  It follows that if this sum is finite, point evaluation at $y$ is a bounded linear functional in $\PP$, i.e.\ 
there is a constant $C_y$
so that $|p(y)|^2 \leq C_y L(p \overline{p})$.  Riesz's theorem provides representation of this functional by the element
$K(x,y)$ as a function in the variable $x$.

Our focus will be on the sum \eqref{E:sum} to show that it converges for every $y \in \C^d$.  In particular we will want
this sum to be uniformly bounded on compact subsets of $\C^d$.  More information on reproducing kernels may be found in
\cite{Sai}.

\section{Main Result}

For $R > 0$ define the $R$-scaling of the multisequence $\{ s_\alpha \}$ to be $\{ \frac{s_\alpha}{R^{|\alpha|}} \}$
and the associated truncated Hankel matrices from this $R$-scaling:
\[
	H_{R,N} = \left( \frac{s_{\alpha + \beta}}{R^{|\alpha + \beta|}} \right)_{0 \leq |\alpha|, |\beta| \leq N}.
\]
Note that if $\{ s_\alpha \}$ is the moment multisequence of a measure $\sigma$, then the $R$-scaling of $s$ is the
moment multisequence of the measure $\sigma_R$, where $\sigma_R(B) = \sigma(RB)$ for every Borel set $B$; we define $RB
= \{ R b : b \in B \}$ in the natural way and define the $R$-scaling of the measure $\sigma$ to be $\sigma_R$.  In the
following theorem, we only assume that $\{s_\alpha\}$ is a positive multisequence, not necessarily that it is a moment
multisequence as we do in the one dimensional case.

We define $\PP_R$ be the Hilbert space associated with the $R$-scaling of $\{ s_\alpha \}$.
If $f(x)$ is in the Hilbert space $\PP$ associated with the sequence $s_\alpha$, then
the mapping $\phi_R : \PP \to \PP_R$ given by $\phi_R(f)(x) = f(Rx)$ is the required isometry.  We will sometimes denote
$\phi_R(f)$ as $f_R$.

Using this isometry, we also relate a system of orthonormal polynomials of the original form to a system of the
scaled form by letting $P_{R,\alpha}(z) = P_\alpha (R z)$.  We choose a set of orthonormal polynomials so that
$\deg(P_\alpha) = \deg(P_{R,\alpha}) = |\alpha|$ (see \cite{DX}).

Let $\lambda_{R,N}$ be the smallest eigenvalue of $H_{R,N}$, and let $\lambda_{R,N} \to \gamma_R$ as $N \to \infty$.  We
now state the main result.

\begin{theorem} \label{T:R}
	Assume that $\lambda_{R,N} > 0$ for all $N \in \N$ and some $R > 0$.  A reproducing kernel for the polynomials exists
	which is uniformly bounded on compact sets if and only if
	$\gamma_R > 0$ for every $R > 0$.
\end{theorem}

\begin{proof}
	This proof follows closely that of Berg, Chen, and Ismail.  The primary difference is in the use of scalings of
	a multisequence.  This is from the fact that it is open problem to whether in the multivariate case the sum
	$\sum_\alpha |P_\alpha(z)|^2$ converging at some nonreal $z$ implies the existence of a reproducing kernel as it does
	in the one variable case.

	We begin by writing the smallest eigenvalue of $H_{R,N}$ as the Rayleigh quotient
	\[
		\lambda_{R,N} = \min \left\{ \sum_{|\alpha| \leq N} \sum_{|\beta| \leq N}
			\frac{s_{\alpha + \beta}}{R^{|\alpha + \beta|}} v_\alpha \overline{v_\beta} : \sum_{|\alpha| \leq N}
			|v_\alpha|^2 = 1 \right\}.
	\]

	Note that if $\lambda_{R,N} > 0$ for some $R > 0$ and $N \in \N$, then $\lambda_{S,N} > 0$ for any $S > 0$.
	We show this by defining $L_R$ to be the functional derived from the multisequence $\left\{
	\frac{s_\alpha}{R^{|\alpha|}} \right\}$, then for any polynomial $p(x)$, $L_R(p(x))$ = $L_S(p(\frac{S x}{R}))$.
	We also write for $p(x) = \sum_{|\alpha| \leq N} a_\alpha x^\alpha$,
	\[
		\lambda_{R,N} = \min \left\{ L_R(|p|^2) : \sum_{|\alpha| \leq N} |a_\alpha|^2 = 1, \: \deg(p) \leq N \right\},
	\]
	therefore if $L_R(|p|^2) = 0$ for some polynomial $p$, then we conclude that $\lambda_{S,N} = 0$ for any $S > 0$.

	We note that the monomials are orthonormal with respect to normalized Lebesgue measure on the unit torus, so we
	rewrite the smallest eigenvalue as
	\[
		\lambda_{R,N} = \min \left\{ L_R(|p|^2) : \int_{\theta \in [0,2 \pi]^d}
			|p(e^{i \theta})|^2 \, \frac{d \theta}{(2 \pi)^d} = 1 \text{, } \deg(p) \leq N \right\}.
	\]
	Since $\lambda_{R,N} > 0$ we write
	\[
		\frac{1}{\lambda_{R,N}} = \max \left\{ \int_{\theta \in [0,2 \pi]^d} |p(e^{i \theta})|^2 \frac{d \theta}{
			(2 \pi)^d} : L_R(|p|^2) = 1 \text{, } \deg(p) \leq N \right\}.
	\]
	We rewrite the polynomial $p(x) = \sum_{|\alpha| \leq N} c_\alpha P_{R,\alpha}(x)$, where $\{P_{R,\alpha}\}$ are
	the standard orthonormal polynomials with respect to a fixed ordering in which $m \geq n \Rightarrow |\alpha_m| \geq
	|\alpha_n|$ for $m,n \in \N$ (see \cite{DX}).  Now considering the matrix $\K_R$, where
	\[
		\K_{R,\alpha,\beta} = \int_{\theta \in [0,2 \pi]^d} P_{R,\alpha} (e^{i \theta}) \overline{P_{R,\beta}(e^{i
			\theta})} \frac{d \theta}{(2 \pi)^d},
	\]
	we express the equation above as another eigenvalue problem:
	\[
		\frac{1}{\lambda_{R,N}} = \max \left\{ \sum_{|\alpha| \leq N} \sum_{|\beta| \leq N} \K_{R,\alpha,\beta}
			c_\alpha \overline{c_\beta} : \sum_{|\alpha| \leq N} |c_\alpha|^2 = 1 \right\}.
	\]
	We also let $\K_{R,N} = (\K_{\alpha,\beta})_{|\alpha|,|\beta| \leq N}$ be a truncation of the matrix $\K_R$.  Note
	that since
	\[
		\int_{\theta \in [0,2 \pi]^d} \left| \sum_{|\alpha| \leq N} c_\alpha P_{R,\alpha}(e^{i \theta}) \right|^2
			\frac{d \theta}{(2 \pi)^d} = \sum_{|\alpha|,|\beta| \leq N}
			c_\alpha \overline{c_\beta} \K_{\alpha,\beta} > 0
	\]
	when $p(z) = \sum_{|\alpha| \leq N} c_\alpha P_{R,\alpha}(z)$ is not the zero polynomial, the matrix $\K_{R,N}$ is
	positive definite.

	($\Rightarrow$) Suppose a reproducing kernel exists, i.e.\
	\[
		\sum_{\alpha \in \N_0^d} |P_\alpha (z)|^2
	\]
	is uniformly bounded on compact subsets of $\C^d$.
	In particular, for $z$ on the torus of radius $R$, there is some
	$M < \infty$ such that $\sum_{\alpha \in \N_0^d} |P_\alpha (z)|^2 \leq M$.  If we consider the $R$-scaling of the
	sequence and the resulting orthonormal polynomials, this implies that
	\[
		\sum_{\alpha \in \N_0^d} |P_{R,\alpha} (z)|^2 \leq M
	\]
	for all $z$ on the unit torus.  Note that the $N^\text{th}$ partial sum is equal to
	the trace of $\K_{R,N}$.  Since $\K_{R,N}$ is positive, the trace is greater than the largest eigenvalue.  Thus
	\begin{align*}
		\frac{1}{\lambda_{R,N}} &\leq \tr(K_{R,N})\\
			&= \sum_{|\alpha| \leq N} \left( \int_{\theta \in [0,2 \pi]^d} |P_{R,\alpha} (e^{i \theta})|^2
				\frac{d \theta}{(2 \pi)^d} \right) \\
			&= \int_{\theta \in [0, 2 \pi]^d} \left( \sum_{|\alpha| \leq N} |P_{R,\alpha} (e^{i \theta})|^2
				\right) \frac{d \theta}{(2 \pi)^d}\\
			&\leq \int_{\theta \in [0,2 \pi]^d} M \frac{d \theta}{(2 \pi)^d}\\
			&= M < \infty.
	\end{align*}
	Hence $\lambda_{R,N} \geq \frac{1}{M}$ for all $N$, and thus is bounded away from zero as $N \to \infty$.

	($\Leftarrow$) Suppose $\lambda_{R,N} \to \gamma_R > 0$ as $N \to \infty$.  From this we conclude that the largest
	eigenvalue of $\K_{R,N}$ is bounded by $\frac{1}{\gamma_R}$ for all $N$, and since it is positive, $\| \K_{R,N} \|$
	is also
	bounded by $\frac{1}{\gamma_R}$.  We wish to show that for each compact $E \subseteq \C^d$, there is some $M <
	\infty$ so that $\sum_{\alpha \in \N_0^d} |P_\alpha(z)|^2 < M$ for all $z \in E$.

	For an arbitary set of complex numbers $\{c_\alpha\}_{|\alpha| \leq N}$, the boundedness of the operator $\K_{R,N}$
	implies
	\[
		\sum_{|\alpha|,|\beta| \leq N} c_\alpha \overline{c_\beta} \K_{R,\alpha,\beta} \leq \frac{1}{\gamma_R}
			\sum_{|\alpha| \leq N} |c_\alpha|^2,
	\]
	and by letting $p(z) = \sum_{|\alpha| \leq N} c_\alpha P_{R,\alpha}(z)$, this is equivalent to
	\begin{equation}\label{E:BCI}
		\int_{\theta \in [0,2 \pi]^d} |p(e^{i \theta})|^2 \frac{d \theta}{(2 \pi)^d} \leq \frac{1}{\gamma_R}
			L_R(|p|^2).
	\end{equation}

	Now let $E$ be a compact set in $\C^d$ and let $R > 0$ be a real number such that $E$ is a subset of the open
	polydisk $D(0,R)^d$.  Then the set $E_R := \frac{1}{R} E$ is contained in the unit polydisk and if $p$ is a
	polynomial, then
	\[
		p(z) = p_R \left( \frac{z}{R} \right),
	\]
	so values taken by $p(z)$ on $E$ are the same as values taken by $p_R(z)$ on $E_R$.

	Let $y \in E_R$.  By Cauchy's integral formula,
	\[
		p_R(y) = \frac{1}{(2 \pi)^d} \int_{[0,2 \pi]^d} \frac{p_R(e^{i \theta})}{(y_1 - e^{i \theta_1})
			(y_2 - e^{i \theta_2}) \cdots (y_d - e^{i \theta_d})} e^{i (\theta_1 + \theta_2 + \cdots + \theta_d)} \,
			d\theta.
	\]
	Using H\"older's inequality, we obtain
	\[
		|p_R(y)|^2 \leq \int_{[\theta \in 2 \pi]^d} |p_R(e^{i \theta})|^2 \frac{d \theta}{(2 \pi)^d} \cdot \int_{\theta
			\in [0,2 \pi]^d} \frac{1}{|y_1 - e^{i \theta_1}|^2 \cdots |y_d - e^{i \theta_d}|^2} \frac{d \theta}{
			(2 \pi)^d}.
	\]
	Letting
	\[
		M = \max \left\{ \frac{1}{\gamma_R} \int_{\theta \in [0,2 \pi]^d} \frac{1}{|w_1 - e^{i \theta_1}|^2 \cdots
		|w_d - e^{i \theta_d}|^2} \frac{d \theta}{(2 \pi)^d} : w \in E_R \right\},
	\]
	we combine this with~\eqref{E:BCI}, so for any polynomial $p_R$,
	\[
		|p_R(y)|^2 \leq M L_R( |p_R(x)|^2 ).
	\]

	Now we pick the particular polynomial $p_R(z) = \sum_{|\alpha| \leq N} \overline{P_{R,\alpha}(y)} P_{R,\alpha}(z)$
	and apply this inequality:
	\[
		\left| \sum_{|\alpha| \leq N} |P_{R,\alpha}(y)|^2 \right|^2 \leq M \sum_{|\alpha| \leq N} |P_{R,\alpha}(y)|^2,
	\]
	Which becomes
	\[
		\sum_{|\alpha| \leq N} |P_{R,\alpha}(y)|^2 \leq M.
	\]
	Since there is a uniform bound over all $N$, this implies
	\[
		\sum_{\alpha \in \N_0^d} |P_{R,\alpha}(y)|^2 \leq M,
	\]
	thus the sum $\sum_\alpha |P_{R,\alpha}(z)|^2$ is bounded uniformly in $\frac{1}{R} E$ which implies that the sum
	$\sum_\alpha |P_\alpha(z)|^2$ is bounded uniformly by $M$ in $E$, thus such a reproducing kernel exists.
\end{proof}

This theorem along with some considerations from the one variable case gives the proof of Theorem \ref{T:BCI}.  In this
case $\{s_n\}$ is a positive multisequence exactly when there is a positive measure $\mu$ which represents $L$ in the
sense that
\[
	L(p) = \int_\R p(x) \, d\mu.
\]
In one dimension it is known that the sum $\sum_{n = 0}^\infty |P_n(z)|^2 < \infty$ for some $z \in \C \backslash \R$ if
and only if the sum converges uniformly on compact subsets of $\C$.  Also in this case, this sum converges if and only if
the measure $\mu$ associated with the moment sequence is indeterminate (see \cite{A}).

So assuming that $\mu$ on $\R$ is indeterminate, then there exists a reproducing kernel for the Hilbert space of
polynomials implying that $\sum_{n = 0}^\infty |P_n(z)|^2$ converges uniformly on compact subsets of $\C$.  By Theorem
\ref{T:R} this leads us to the conclusion that $\gamma_1 > 0$, i.e.\ the smallest eigenvalues of $H_N$ are uniformly
bounded away from 0.

Now assuming that $\gamma_1 > 0$, the proof of Theorem \ref{T:R} implies that the sum $\sum_{n = 0}^\infty |P_n(z)|^2$
is bounded uniformly on
compact subsets of the open unit disk.  Thus there is some $z_0 \in \mathbb{D} \backslash \R$ so that $\sum_{n = 0}^\infty
|P_n(z_0)|^2 < \infty$, hence the measure $\mu$ is indeterminate.

\section{Application}

We would like to know more about $\gamma_R$.  If we consider it as a function $\gamma : (0,\infty) \to [0,\infty)$ given by
$\gamma(R) = \gamma_R$, then Theorem \ref{T:R} implies that if $\gamma(R)> 0$, then $\gamma(S) > 0$ for all $S < R$.  We
see this since $\gamma(R) > 0$ implies that a reproducing kernel for $\PP$ exists and converges uniformly on compact
subsets of the polydisk of radius $R$.  Then in the converse portion of the argument, the fact that a reproducing kernel
exists which converges uniformly on the torus of radius $S$ is used to show that $\lambda_{S,N}$ are uniformly bounded away
from zero as $N \to \infty$, hence $\gamma(S) > 0$.

If a multisequence in any number of variables can be represented by a measure $\sigma$, then $\sigma$ is indeterminate if
and only if every $R$-scaling $\sigma_R$ is as well.  In one
variable, Berg, Chen, and Ismail's result implies that $\gamma_1 > 0 \Leftrightarrow \gamma_R > 0$ for some $R > 0$, which
is equivalent to $\gamma_S > 0$ for every $S > 0$.  This is an open question whether this property holds for multisequences
in multiple variables, but in Theorem \ref{T:R2} we establish a class of multisequences which do possess this property.

We would like to extend more of Theorem \ref{T:BCI} to the $d > 1$ case.  Presently we have the conditions
\begin{align}
	&\gamma_1 > 0, \label{c1}\\
	&\gamma_R > 0, \text{ for all $R > 0$}, \label{c2}\\
	&\sum_{\alpha \in \N_0^d} |P_\alpha (z)|^2 \text{ is bounded uniformly on compact sets}, \label{c3}\\
	&\sum_{\alpha \in \N_0^d} |P_\alpha (z)|^2 < \infty \text{ for some $z \in (\C \backslash \R)^d$}, \label{c4}
\end{align}
are related by the implications
\[
	\eqref{c1} \Leftarrow \eqref{c2} \Leftrightarrow \eqref{c3} \Rightarrow \eqref{c4},
\]
and also $\eqref{c1} \Rightarrow \eqref{c4}$ if we pick some $z$ contained in the open unit polydisk.  We would like also
to have $\eqref{c1} \Rightarrow \eqref{c2}$ and
$\eqref{c4} \Rightarrow \eqref{c3}$ in $d > 1$ as it is in $d = 1$.  We show this is the case for a certain class of
multisequences.

\begin{theorem} \label{T:R2}
	Let $\{s_\alpha\}$ be a positive multisequence which satisfies the condition
	\begin{equation} \label{mult}
		s_\alpha = s_{(\alpha_1, 0, 0, \ldots, 0)} s_{(0,\alpha_2, 0, \ldots, 0)} \cdots s_{(0, \ldots, 0, \alpha_d)}.
	\end{equation}
	Then there is a positive measure $\mu$ on $\R^d$ which represents the
	multisequence in the sense that
	\[
		\int_{\R^d} x^\alpha \, d\mu = s_\alpha;
	\]
	the conditions \eqref{c1}, \eqref{c2}, \eqref{c3}, and \eqref{c4} are equivalent; and these conditions
	all imply indeterminacy of $\mu$.
\end{theorem}

\begin{proof}
	Let $\{ s_\alpha \}$ satisfy this condition and let $L$ be the functional derived from the multisequence.  We define
	a collection of $d$ sequences by for each $j = 1, \ldots, d$, let
	\[
		s_{j,n} = s_{(0, \ldots, 0, n, 0, \ldots, 0)},
	\]
	where the $n$ is in the $j^\text{th}$ place.  The positivity of the sequence $s_{j,n}$ follows from the positivity of
	the multisequence $\{ s_\alpha \}$, so for each $1 \leq j \leq d$, the theorem of Hamburger (\cite{A} p.30) implies
	that there exists a positive measure $\mu_j$ on $\R$ for which
	\[
		\int_\R x^n \, d\mu_j(x) = s_{j,n}.
	\]

	If we define a measure $\mu := \mu_1 \otimes \mu_2 \otimes \cdots \otimes \mu_d$ on $\R^d$, the multiplicative
	condition \eqref{mult} implies the first assertion, that
	\[
		\int_{\R^d} x^\alpha \, d\mu = \left( \int_\R x_1^{\alpha_1} \, d\mu_1 \right) \cdots \left( \int_\R
			x_d^{\alpha_d} \, d\mu_d \right) = s_{1,\alpha_1} s_{2,\alpha_2} \cdots s_{d,\alpha_d} = s_\alpha,
	\]
	for each monomial $x^\alpha$.

	For the measure $\mu$, we construct a suitable set of orthonormal polynomials $\{ P_\alpha (x)\}$.  For each
	$1 \leq j \leq d$, let
	$\{ P_{j,n}(x) \}_{n = 0}^\infty$ be the standard set of orthonormal polynomials associated with the measure $\mu_j$
	in $\R$ (\cite{A} p. 3).  If we define for each $\alpha \in \N_0^d$ the polynomial
	\[
		P_\alpha (x) := P_{1,\alpha_1}(x_1) P_{2,\alpha_2}(x_2) \cdots P_{d, \alpha_d}(x_d),
	\]
	then the resulting set of polynomials $\{ P_\alpha (x) \}$ span $\C[x]$ and satisfy the orthonormal property
	\[
		\int_{\R^d} P_\alpha(x) \overline{P_\beta(x)} \, d\mu(x) = \delta_{\alpha_1, \beta_1} \delta_{\alpha_2,
			\beta_2} \cdots \delta_{\alpha_d, \beta_d} = \delta_{\alpha, \beta}.
	\]

	Now for any $z \in \C^d$ and $N \in \N$, we consider the sum
	\begin{align} \label{mult2}
		\sum_{\overset{|\alpha_j| \leq N}{1 \leq j \leq d}} |P_\alpha(z)|^2 = \left( \sum_{\alpha_1 = 0}^N
			|P_{1,\alpha_1}(z_1)|^2 \right)
			&\left( \sum_{\alpha_2 = 0}^N |P_{2,\alpha_2}(z_2)|^2 \right) \cdots \notag \\
			&\cdots \left( \sum_{\alpha_d = 0}^N |P_{d,\alpha_d}(z_d)|^2 \right).
	\end{align}
	This sum converges as $N \to \infty$ if and only if each of the sums on the right hand side converges.  We will show
	that for this measure, $\eqref{c1} \Rightarrow \eqref{c4}$ for a generic $z \in (\C \backslash \R)^d$ and $\eqref{c4}
	\Rightarrow \eqref{c3}$.  In the process we prove that these imply indeterminacy of $\mu$.

	($\eqref{c1} \Rightarrow \eqref{c4}$).  Suppose $\gamma_1 > 0$.  Define
	\[
		J_N^j = \left( s_{(0, \ldots, 0, k + l, 0, \ldots, 0)} \right)_{k,l = 0}^N,
	\]
	where the $k + l$ is in the $j^\text{th}$ coordinate, to be the $(N + 1)
	\times (N + 1)$ truncated Hankel matrix associated with the measure $\mu_j$, and let $\eta_N^j$ be its smallest
	eigenvalue.  Then $\eta_N^j > \gamma_1$ for each $j$ and $N \in \N$ since $J_N^j$ is a compression of the operator
	$H_N$.  By Theorem \ref{T:BCI} it follows that $\mu_j$ is indeterminate and
	\[
		\sum_{\alpha_j = 0}^\infty |P_{j,\alpha_j}(z_j)|^2
	\]
	converges for some $z_j \in \C \backslash \R$.  If $z = (z_1, z_2, \ldots, z_d) \in (\C \backslash \R)^d$, then
	\eqref{c4} follows from \eqref{mult2}.

	($\eqref{c4} \Rightarrow \eqref{c3}$).  Now let $z \in (\C \backslash \R)^d$ be a point such that $\sum_{\alpha \in
	\N_0^d} |P_\alpha (z)|^2 < \infty$.  Via \eqref{mult2}, this implies that $\sum_{\alpha_j = 0}^\infty
	|P_{j,\alpha_j}(z_j)|^2 < \infty$ for each $1 \leq j \leq d$, which implies that $\mu_j$ is indeterminate.  Then
	$\sum_{\alpha_j = 0}^\infty |P_{j,\alpha_j} (z_j)|^2$ converges uniformly on compact subsets of $\C$.  Let $E$ be
	a compact subset of $\C^d$ and $L \geq 0$ so that $E \subseteq \left( \overline{\mathbb{D}(0,L)} \right)^d$.  Then
	\eqref{mult2} tells us that $\sum_{\alpha \in \N_0^d} |P_\alpha (z)|^2$ converges uniformly on $\left(
	\overline{\mathbb{D}(0,L)} \right)^d$ and thus on $E$.

	We now have that \eqref{c1}, \eqref{c2}, \eqref{c3}, and \eqref{c4} are equivalent when \eqref{mult} applies, and our
	arguments have shown that these conditions imply indeterminacy for each $\mu_j$.  Since at least one of these
	measures is indeterminate, Petersen's theorem leads to indeterminacy of $\mu$.
\end{proof}

This theorem provides that when a positive multisequence is multiplicative in this sense, a solution to the associated
moment problem exists, and the existence of a reproducing kernel implies indeterminacy.  In the case of deciding
indeterminacy of an arbitrary moment multisequence, the role of reproducing kernels is not yet fully clear.  Fuglede's
introduction of the notions of strong determinacy and ultradeterminacy \cite{F} and the creation of new indeterminate
moment problems from old ones \cite{PuSc} provide ground for further study, and we will address these in a future article.

\bibliographystyle{amsplain}

\end{document}